\documentclass{amsart}
\usepackage{amssymb,amsfonts,amsmath}
\usepackage{graphics}
\usepackage{graphicx}
\usepackage{epsfig}
\usepackage{color}
\usepackage{enumerate}
\usepackage{amsthm}

\newtheorem{lema}{Lemma}[section]

\newtheorem*{Domar teo}{Domar's Theorem}

\newtheorem{proposicion}[lema]{Proposition}

  {\par \hfill \fbox{}}
  {\par \hfill \fbox{}}

\def\CC{{\mathbb C}}
\def\D{{\mathcal D}}
\def\DD{{\mathbb D}}

\def\RR{{\mathbb R}}
\def\TT{{\mathbb T}}

\def\A{\mathcal A}
\def\C{\mathcal C}
\def\D{\mathcal D}
\def\H{\mathcal H}
\def\M{\mathcal M}
\def\MD{\mathcal{M}(\mathcal{D})}
\def\Mphi{\mathcal{M}(\phi)}
\def\lesim{\lesssim }

\def\epsilon{\varepsilon}

\def\phi{\varphi}

\newtheorem{thm}{Theorem}[section]

\newtheorem{exam}{Example}[section]

\def\beginpf{\begin{proof}}
\def\endpf{\end{proof}}

\def\beq{\begin{equation}}
\def\eeq{\end{equation}}

  {\par \hfill }

\begin{document}

\title[Weighted composition operators]{Weighted composition operators on the Dirichlet space: boundedness and spectral properties}
\author{I. Chalendar}
\address{Universit\'e de Lyon; CNRS;\newline
Universit\'e Lyon 1; INSA de Lyon;
Ecole Centrale de Lyon, \newline
CNRS, UMR 5208, Institut Camille Jordan\newline
43 bld. du 11 novembre 1918, F-69622 \newline
Villeurbanne Cedex, France}
\email{chalendar@math.univ-lyon1.fr}
\author{E. A. Gallardo-Guti\'errez}
\address{
Universidad Complutense de Madrid e ICMAT\newline
Departamento de An\'alisis Matem\'atico,\,\newline
Facultad de Ciencias Matem\'aticas,\newline
Plaza de Ciencias 3\newline
28040, Madrid (SPAIN)}
\email{eva.gallardo@mat.ucm.es}
\author{J. R. Partington}
\address{School of Mathematics,\newline
 University of Leeds,\newline
Leeds LS2 9JT, U.K.} \email{J.R.Partington@leeds.ac.uk}

\thanks{The authors are partially supported by Plan Nacional I+D grant no.
MTM2013-42105-P}

\subjclass{Primary 47B38.}
\keywords{Dirichlet space, weighted composition operator}
\date{September 2014, Revised: February 2015}


\begin{abstract}
Boundedness of weighted composition operators
$W_{u,\varphi}$ acting on the classical Dirichlet space
$\mathcal{D}$ as $W_{h,\varphi}f= h\, (f\circ \varphi)$ is studied in terms of the
\emph{multiplier space} associated to the symbol $\varphi$, i.e.,
$\Mphi=\{ u \in \D: W_{u,\phi} \hbox{ is bounded on } \D \}$. A prominent role is played by the multipliers of
the Dirichlet space. As a consequence, the spectrum of $W_{u,\varphi}$ in $\mathcal{D}$
whenever $\varphi$ is an automorphism of the unit disc is studied, extending a recent work of
Hyv\"arinen, Lindstr\"om,  Nieminen and  Saukko \cite{HLNS} to the context of the Dirichlet space.
\end{abstract}

\maketitle

\section{Introduction and Preliminaries}

Let $\DD$ denote the open unit disc in the complex plane $\CC$. The Dirichlet space $\D$ consists of analytic functions $f$ on $\DD$
such that the norm
$$ \|f\|_{\mathcal{D}}^2=|f(0)|^2 + \int_{\mathbb{D}} |f'(z)|^2 \,
dA(z) $$
\noindent is finite. Here $A$ stands for the normalized Lebesgue
area measure of the unit disc. Observe that for a univalent
function $f$, the integral above is just the area of
$f(\mathbb{D})$.

It is well known that $\D \subset \H^2 \subset \A^2$, where $\H^2$ and $\A^2$ denote
respectively the Hardy and Bergman spaces on $\DD$, and that $f \in \D$ if and only
if $f' \in \A^2$. The recent monograph \cite{EKMR} is an excellent source
to learn about the Dirichlet space and its particular issues.

If $\varphi$ is an analytic function on $\mathbb{D}$ with
$\varphi(\mathbb{D})\subset \mathbb{D}$, then the equation $$
C_{\varphi} f =f\circ \varphi $$ defines a composition operator
$C_{\varphi}$ on the space of all holomorphic functions on the
unit disc $\mathcal{H}(\mathbb{D})$.  On the
Dirichlet space $\mathcal{D}$, a necessary condition for
$C_{\varphi}$ to be bounded is that $\varphi\in \mathcal{D}$. Nevertheless,
not all the Dirichlet functions induce bounded composition operators on $\mathcal{D}$. Such
functions were characterized in 1980 by C. Voas \cite{Vo} in his Ph.D.\ thesis.

In this work, we shall be concerned with weighted composition operators on $\D$: for $u \in \D$ and $\phi$ a holomorphic self-map of $\DD$ we define the weighted composition operator $W_{u,\phi}$ on $\D$ by
\[
(W_{u,\phi}f)(z)=u(z) f(\phi(z)),
\]
noting that $W_{u,\phi}$ is not, in principle, a bounded operator on $\D$. It is clear that if $C_{\varphi}$
is a bounded operator on $\D$ and $u$ is a \emph{multiplier} of $\D$, that is,
the Toeplitz operator $T_u: f \mapsto u   f$ is defined everywhere on $\D$ and hence bounded, the
weighted composition operator $W_{u,\phi}$ on $\D$ is obviously bounded.

A well known fact about the  Dirichlet space is that the algebra $\MD$ consisting of the multipliers of $\D$ is not that easy to describe.
Indeed, their elements were characterized by Stegenga \cite{St} in a remarkable paper in terms of a condition involving the logarithmic capacity of their boundary values. In particular, the strict inclusion $\MD \subset \D \cap \H^\infty$ holds. Here $\H^{\infty}$ denotes the space
of bounded analytic functions in $\DD$ endowed with the sup-norm. A straightforward reformulation in terms of Carleson measures for $\D$ (that is, there is a continuous injection from $\D$ into $L^2(\DD,\mu)$), yields the fact  that  $u \in \MD$ if and
only if $u$ is bounded and the measure $\mu$ defined by $d\mu(z)=
|u'(z)|^2 \, dA(z)$ is a  Carleson measure for $\D$. We refer to \cite{Wu} for multipliers and Carleson measures in Dirichlet spaces (and to \cite{alsaker,EKMR} for more on the subject of multipliers).

Concerning boundedness of weighted composition operators on $\D$, let us remark that one may construct self-maps of the unit disc $\varphi$ such that $\varphi\not \in \D$ and a multiplier $u\in \MD$
such that $W_{u,\varphi}$ is bounded in the Dirichlet space. For instance, let
$u(z)=(1-z)^2$ and let $\phi$ be the infinite Blaschke product with zeroes $(1-1/n^2)_{n \ge 1}$. Now $\phi\not\in \D$, so $C_\phi$
is clearly unbounded.  However, for $f \in \D$ we have
\[
((1-z)^2(f \circ \phi))'= -2(1-z)(f \circ \phi) + (1-z)^2 (f' \circ \phi) \phi';
\]
now the first term is clearly in the Bergman space $\A^2$, while for the
second term we observe that
$(1-z)^2 \phi'$ is a bounded analytic function in $\DD$ and $(f' \circ \phi) \in \A^2$, so that it also lies in $\A^2$ (with control of norms),
showing that $W_{u,\phi}$ is bounded on $\D$.

Therefore, facing the problem of describing the weighted composition operators taking $\D$ boundedly into itself deals not only with
the multipliers of $\D$ but also with those self-maps of the unit disc that may induce unbounded composition operators in $\D$.

At this regards, for a self-map $\phi$ of the unit disc $\DD$, we define the \emph{multiplier space} $\Mphi$ associated to $\phi$ by
\[
\Mphi=\{ u \in \D: W_{u,\phi} \hbox{ is bounded on } \D \}.
\]
It is clear that if $C_\phi$ is bounded on $\D$, then  $\MD \subseteq \Mphi \subseteq \D$. Moreover, if $\phi$ induces an unbounded $C_{\phi}$ in $\D$, then $\MD$ is no longer contained in $\Mphi$ since, in such a case, this latter space does not contain the constant functions.

The aim of this work is twofold. On one hand, we are interested in identifying the multiplier space $\Mphi$ for self-maps $\varphi$ of $\DD$.
Indeed, we will be able to characterize the extreme cases whenever $\phi$ is a self-map of $\DD$ belonging to $\D$. Let us remark here that in the case of the Hardy space $\H^2$, where $C_\phi$ is automatically bounded and the multiplier space is $\M (\H^2)=\H^\infty$, Gallardo-Guti\'errez, Kumar and Partington proved that $\M_{\H^2}(\phi)=\H^2$ if and only if $\|\phi\|_\infty<1$ (see \cite{GKP}) and $\M_{\H^2}(\phi)=\H^\infty$ if and only if $\phi$ is a finite Blaschke product; this latter statement was showed previously in a different way in \cite{CH} and \cite{Matache}.

Let us also point out that in the course of our findings, we will prove a \emph{Decomposition Theorem} for the Dirichlet space (cf. Theorem \ref{thm:decomp}), which is interesting in its own and whose proof is based on the theory of model spaces for the shift operator in the Hardy space (see \cite{Ni} for more information about model spaces). As far as we know, this is the first time model spaces come into play with the Dirichlet space.

On the other hand, we are interested in the spectral properties of weighted composition operators in $\D$. In \cite{higdon}, Higdon computed the spectrum of composition operators in $\D$ induced by linear fractional self-maps of $\DD$. The techniques developed there were quite different from those carried over by Cowen in \cite{Co83} in the corresponding case of the Hardy space $\H^2$ (see also \cite[Chapter 7]{CMc}), due to the particular nature of $\D$.

In a very recent work,  Hyv\"arinen,  Lindstr\"om, Nieminen and Saukko \cite{HLNS} have described the spectra  of invertible weighted composition operators $W_{u,\phi}$ acting on a large class of analytic function spaces including the weighted Bergman and the weighted Hardy spaces; generalizing previous results obtained in \cite{Gu}. Nevertheless, as they also remark, their results do not apply directly to the Dirichlet space since they rely on the fact that the
algebra of the multipliers of the spaces considered is $\H^\infty$. Our aim is to extend Hyv\"arinen,  Lindstr\"om, Nieminen and Saukko's results to the context of the Dirichlet space $\D$, pointing out that their techniques are no longer working in  $\D$.

The last section of the paper gives a description of the spectra of invertible weighted composition operators. We first note (see Proposition \ref{prop:31}) that a bounded weighted composition operator $W_{u,\varphi}$  in the Dirichlet space $\D$ is invertible if and only if $u$ is a multiplier bounded away from zero in $\DD$ and $\varphi$ is an automorphism of the unit disc. Consequently, three separate cases are considered, depending on the nature of the disc automorphism $\varphi$: elliptic, parabolic or hyperbolic. When $\varphi$ is parabolic, causal operators will play a prominent role in order to determine explicitly the spectrum of $W_{u,\varphi}$.

\section{Boundedness of weighted composition operators}

In this section, we study boundedness of weighted composition operators in the Dirichlet space. In order to show the results at this respect,
we prove a \emph{Decomposition Theorem} for $\D$ based on model spaces.

Let $B$ be a finite Blaschke product  and write $K_B$ for the model space $K_B=\H^2 \ominus B\H^2$,
which is finite-dimensional; indeed $\dim K_B=\deg B$. Observe that if $g_k \in K_B$, then it does not matter
which norm we use, since $K_B$ is finite-dimensional and all norms are equivalent. We proceed to state the \emph{Decomposition Theorem}
in its full generality, since the main arguments of the proof also work for the Bergman  space.

\begin{thm}[\textbf{Decomposition Theorem}]\label{thm:decomp}
Let $B$ be a finite Blaschke product such that $B(0)=0$. Then
\begin{enumerate}
\item $f \in \H^2$ if and only if $f = \sum_{k=0}^\infty g_k B^k$ (convergence in $\H^2$ norm)
with $g_k \in K_B$ and
$\sum_{k=0}^\infty \|g_k\|^2 < \infty$.
\item $f \in \D$ if and only if $f = \sum_{k=0}^\infty g_k B^k$ (convergence in $\D$ norm) with $g_k \in K_B$ and
$\sum_{k=0}^\infty(k+1)\|g_k\|^2 < \infty$.
\item $f \in \A^2$
if and only if $f = \sum_{k=0}^\infty g_k B^k$  (convergence in $\A^2$ norm) with $g_k \in K_B$ and
$\sum_{k=0}^\infty \|g_k\|^2/(k+1) < \infty$.
\end{enumerate}
\end{thm}

Statement (1) is included for the sake of completeness since it is  a standard fact that $\H^2=K_B \oplus BK_B \oplus B^2K_B \oplus...$ as an orthonormal direct sum.

\emph{A word about notation}. Throughout this work, $a\lesssim b$ will denote that there exists an independent constant $C$ such that
$a\leq C\, b$; this constant may be different in each instance.

\beginpf
We proceed to prove (2) and (3). We claim that for finite sums
\[
\| \sum_{k=0}^N g_k B^k \|_\D^2 \lesim \sum_{k=0}^N (k+1) \|g_k\|^2 \quad \hbox{and} \quad
 \| \sum_{k=0}^N g_k B^k \|_{\A^2}^2 \lesim \sum_{k=0}^N \|g_k\|^2/(k+1)
\]
with the implied constants independent of $N$.

Let $e_1,\ldots,e_n$ be a basis of the space $K_B$. Then, writing $g_k=\sum_{\ell=1}^n a_{k\ell} e_\ell$ we have,
since $C_B$ is a bounded composition operator in $\D$, that
\begin{eqnarray*}
\left \| \sum_{k=0}^N \sum_{\ell=1}^n a_{k\ell}e_\ell B^k \right \|_\D & = & \left \| \sum_{\ell=1}^n e_\ell \sum_{k=0}^N  a_{k\ell} B^k\right \|_{\D}\\
& \le& \sum_{\ell=1}^n \|e_\ell\|_{\MD} \left \| \sum_{k=0}^N a_{k\ell} B^k \right \|_\D \\
& \lesim & \sum_{\ell=1}^n \left( \sum_{k=0}^N (k+1) |a_{k\ell}|^2 \right)^{1/2} \\
&\lesim & \left( \sum_{\ell=1}^n  \sum_{k=0}^N (k+1) |a_{k\ell}|^2 \right)^{1/2},
\hbox{ by equivalence of norms on } \CC^n,\\
&\lesim & \left( \sum_{k=0}^N (k+1) \|g_k\|^2 \right)^{1/2},
\end{eqnarray*}
where the notation $\|e_\ell\|_{\MD}$ represents $\sup\{\|e_\ell f\|_\D: \|f\|_\D\leq 1\}$.
A similar calculation can be made in $\A^2$.

To obtain the converse inequality, we use the dual pairing between $\D$, equipped with the equivalent
norm $\| \sum_{k=0}^\infty a_k z^k\|_{X}^2 = \sum_{k=0}^\infty (k+1) |a_k|^2$
and $\A^2$, equipped with the equivalent norm $\| \sum_{k=0}^\infty a_k z^k\|_{Y}^2 = \sum_{k=0}^\infty |a_k|^2/(k+1)$,
namely
\[
\left \langle \sum_{k=0}^\infty a_k z^k,  \sum_{k=0}^\infty b_k z^k \right \rangle =  \sum_{k=0}^\infty a_k \overline{b_k}.
\]

For a finite sum $\sum_{k=0}^N g_k B^k$ we take $h_k \in K_B$ with
$h_k=(k+1)g_k$ for each $k$.
Now we use the $\H^2$ orthogonality of $ B^r K_B $ and $ B^s K_B $ for positive integers $r$ and $s$ with $r \ne s$, to deduce that
\[
\left \langle \sum_{k=0}^N g_k B^k, \sum_{k=0}^N h_k B^k \right \rangle = \sum_{k=0}^N (k+1)\|g_k\|_{\H^2}^2 .
\]
Since
\[
\left \| \sum_{k=0}^N h_k B^k \right\|^2_Y \lesim \sum_{k=0}^N   \|h_k\|^2/(k+1)=\sum_{k=0}^N (k+1)\|g_k\|^2,
\]
it follows that
\[
\left\| \sum_{k=0}^N g_k B^k \right \|^2_X \gtrsim \sum_{k=0}^N (k+1) \|g_k\|^2,
\]
and so we have a uniform equivalence of the Dirichlet norm and the quantity
$$\left(\sum_{k=0}^N (k+1)^2 \|g_k\|^2\right)^{1/2},$$
at least for finite sums.
Since the Dirichlet space is contained in the Hardy space we may make the obvious extension to the
whole of $\D$ using infinite sums.

The argument for the Bergman space is analogous: once more we have an equivalence of norms, and since
the Hardy space is dense in the Bergman space we obtain the required result.
\endpf

Recall that if $C_\phi$ is bounded then $\MD \subseteq \Mphi \subseteq \D$. For $\phi$ a finite Blaschke product the
space of weighted composition operators is as small as possible, as the following result shows.

\begin{thm}\label{Theorem:inner}
Let $\phi$ be an inner function. Then $\M(\phi)=\M(\D)$ if and only if $\phi$ is a finite
Blaschke product.
\end{thm}

\beginpf
If $\phi$ is inner but not a finite Blaschke product, then $\phi \not\in \D$, and so $C_\phi$ is unbounded. Thus,
taking $u(z) \equiv 1$, we have that $u \in \M(\D)$ but $u \not \in \M(\phi)$.

Now suppose that $\phi=B$, a finite Blaschke product. Let $T_u$ denote the map $f\in \D \mapsto u   f$.  We must show that $T_u C_B$ is bounded if and only if $u \in \M(\D)$.
It is clear that if $u \in \M(\D)$ then $T_u C_B$ is bounded on $\D$, since both $T_u$ and $C_B$ are
bounded.

For the converse, note that without loss of generality we may take $B(0)=0$, since if for some $a$ we have $B(a)=0$, then,
setting $$\phi_a(z)=\frac{a-z}{1-\overline a z}$$ we have $B \circ \phi_a(0)=0$ and
$$C_{\phi_a} T_u C_B = T_{u\circ \phi_a}C_{B \circ \phi_a}.$$ Now
$T_u C_B$ is bounded if and only if $T_{u\circ \phi_a}C_{B \circ \phi_a}$ is bounded, and showing that $u \in \M(\D)$ is
equivalent to showing that $u \circ \phi_a \in \M(\D)$.

Now, given that $T_u C_B$ is bounded, let $f \in \D$. By Theorem \ref{thm:decomp} we may write
\[
f=\sum_{k=0}^\infty g_k B^k = \sum_{k=0}^\infty \sum_{\ell=1}^n a_{k\ell}e_\ell B^k,
\]
where each $g_k \in K_B$ , $\{e_1,\ldots,e_n\}$ is a basis of $K_B$, and for each $\ell$ the scalars $(a_{k\ell})$ satisfy
\[
\sum_{k=0}^\infty (k+1)|a_{k\ell}|^2 \lesim \sum_{k=0}^\infty (k+1)\|g_k\|^2 \lesim \|f\|^2_\D.
\]
Write $f_\ell(z)=\sum_{k=0}^\infty�a_{k\ell} z^k$ so that $\|f_\ell\|_\D \lesim \|f\|_\D$.
Then
\[
\left \| u \sum_{k=0}^\infty a_{k\ell} B^k  \right \|_\D = \|T_u C_B f_\ell\|_\D \lesim \|T_u C_B\| \| f\|_\D.\]
 Now $e_1,\ldots,e_n$ lie in $\M(\D)$ and we conclude that
\[
 \|uf\|_\D = \left \|  \sum_{\ell=1}^n e_\ell T_u C_B f_\ell \right \|_\D \lesim \|T_u C_B\| \|f\|_\D,
\]
so that $u \in \M(\D)$ and the Theorem is proved.
\endpf

Next example shows that the assumption about $\varphi$ being inner cannot be relaxed; even if $\| \phi \|_\infty=1$ and $C_{\varphi}$ is bounded in $\D$.

\begin{exam}{\rm
We can have $\M(\phi) \ne \M(\D)$ even if $\| \phi \|_\infty=1$ and $C_{\varphi}$ is bounded in $\D$. Let us consider
\[
\phi(z)=\frac{1-z}{2} \qquad \hbox{and} \qquad u(z)=\sum_{k=2}^\infty \frac{z^k}{k(\log k)^{3/4}}.
\]
Note that $u \in \D \setminus  \M(\D)$ (see \cite[Thm.~5.1.6]{EKMR}). Nonetheless, $T_u C_\phi$ is
bounded; that is, $u \in \M(\phi)$.

In order to show that, let $f \in \D$. We have
$$(u(f \circ \phi))'=u'(f \circ \phi) + u(f' \circ \phi)\phi'.$$
It suffices if we prove that each of these terms lies in the Bergman space $\A^2$.



Let us split the disc into disjoint measurable sets $\DD=D_1 \cup D_2$, where $D_1$ is a small neighbourhood of $1$,
mapped by $\phi$ into a small disc about $0$.

	\includegraphics[width=20cm, trim = 5mm 220mm 0 5mm, clip  ]{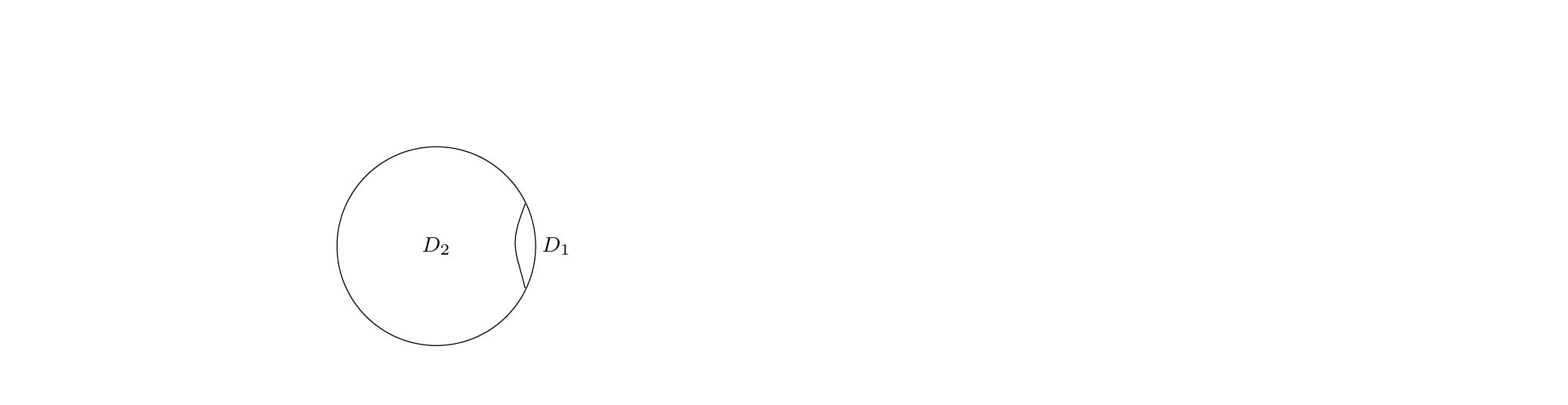}
\begin{center}
A partition of the unit disc.
\end{center}

\medskip

On $D_1$ both $u$ and $u'$ are square-integrable with respect to Lebesgue measure, while $f \circ \phi$ and
$f' \circ \phi$ are uniformly bounded by constants depending only on
the norm of $f$. Likewise, on $D_2$, $u$ and $u'$ are bounded, while $f \circ \phi$ and $f' \circ \phi$ are square-integrable.
Our conclusion is that $T_u C_\phi$ is bounded on the Dirichlet space $\D$.

}
\end{exam}

At the other extreme, we have the following result.

\begin{thm}
Let $\phi$ be an analytic self-map of $\DD$. Then $\M(\phi)=\D$ if and only if
\begin{enumerate}
\item $\|\phi\|_\infty < 1$, and
\item $\phi \in \M(\D)$.
\end{enumerate}
\end{thm}

\beginpf
Suppose first that $\phi$ satisfies conditions (1) and (2). Let $u$ be a Dirichlet function.
In order to prove that $T_u C_\phi$ is a bounded operator on $\D$, let $f\in \D$ and consider
\[
(u( f \circ \phi))'= u'(f \circ \phi) + u (f'\circ \phi)\phi'.
\]
It suffices to show that each of these terms lies in the Bergman space $\A^2$.

We have that $\|f \circ \phi\|_\infty \lesim \|f\|_\D$ since, for $w \in \DD$,
$$ ( f\circ\phi)(w)=\langle f, k_{\phi(w)} \rangle_{\D},$$ where
$k_{\phi(w)}$ denotes the reproducing kernel at $\phi(w)$ in $\D$, which
is bounded in norm independently of $w$ since $\|\phi\|_\infty<1$. Hence, since $u \in \D$, we have
\begin{equation}\label{eq1}
\|u'(f \circ \phi)\|_{\A^2} \lesim \|f\|_\D.
\end{equation}

Also, $\|f' \circ \phi\|_\infty \lesim \|f\|_\D$, since the derivative kernels $ k'_{\phi(w)}$
satisfying $$f'(\phi(w))=\langle f, k'_{\phi(w)} \rangle_{\D}$$
are also uniformly bounded in  norm
when  $\|\phi\|_\infty<1$. Now by means of condition (2), $\phi \in \M(\D)$, and therefore the measure $|\phi'(z)|^2 \, dA(z)$ is a Carleson measure for $\D$. Thus, it follows that
\begin{equation}\label{eq2}
\|u (f'\circ \phi)\phi'\|_{\A^2} \lesim \|f\|_\D.
\end{equation}
From (\ref{eq1}) and (\ref{eq2}), one gets that $T_u C_\phi$ is a bounded operator for all $u \in \D$.

Conversely, if $T_u C_\phi$ is bounded for all $u \in \D$, then by the Closed Graph Theorem $C_\phi$
maps $\D$ boundedly into $\M(\D)$ and hence into $\H^\infty$. Suppose that $\| \phi\|_\infty=1$; then we may find an unbounded
$f \in \D$, and by considering functions $f_n(z)=f(e^{i\theta_n}z)$ for suitable angles $\theta_n$, obtain a sequence of
normalized functions $f_n$ in $\D$ and $(z_n)\subset \DD$ with $|z_n| \to 1$ such that $|f_n(\phi(z_n))| \to \infty$. This is
a contradiction, so we conclude that $\|\phi\|_\infty<1$.

We now see that, for $f \in \D$ fixed, we have $\|u(f' \circ \phi)\phi'\|_{\A^2} \lesim \|u\|_\D$. So let $f(z)=z$. We
conclude that $T_{\phi'}: \D \to \A^2$ is bounded, or equivalently $|\phi'(z)|^2 dA(z)$ is a Carleson measure for
$\D$. Given that $\phi$ is bounded, this condition implies that $\phi \in \M(\D)$, as mentioned in the introduction. This concludes the proof of the Theorem.
\endpf

\section{Spectral properties}

In this section, we are interested in describing the spectra  of invertible weighted composition operators in the Dirichlet space.
As we pointed out in the introduction, the techniques in \cite{HLNS} depends strongly on the fact that the algebra of the multipliers contains
$\H^{\infty}$.

The next result identifies invertible weighted composition operators in the Dirichlet space.
It can be found, for example, in   \cite[Thm~3.3]{bourdon} and \cite[Cor.~11]{zhao}.

\begin{proposicion}\label{prop:31}
Let $W_{h,\varphi}$ be a bounded weighted composition operator in the Dirichlet space $\D$. Then $W_{h,\varphi}$ is invertible in $\D$ if and only if
$h\in \MD$, bounded away from zero in $\DD$ and $\varphi$ is an automorphism of $\DD$. In such a case, the inverse operator of $W_{h,\varphi}:\D\to \D$ is also a weighted composition operator and
$$
(W_{h,\varphi})^{-1}=\frac{1}{h\circ \varphi^{-1}} C_{\varphi^{-1}}.
$$
\end{proposicion}

Recall that an automorphism $\varphi$ of
$\mathbb{D}$ can be expressed in the form
$$
\varphi(z)=e^{i\theta}\, \frac{p-z}{1-\overline{p}z} \qquad (z\in
\mathbb{D}),
$$
where $p\in \mathbb{D}$ and $-\pi<\theta\leq \pi$. Recall that
$\varphi$ is called  \emph{hyperbolic} if $|p|>\cos(\theta/2)$
(thus, $\varphi$ fixes two points on $\partial \mathbb{D}$);
\emph{parabolic} if $|p|=\cos(\theta/2)$ (so, $\varphi$ fixes
just one point, located on $\partial \mathbb{D}$) and
\emph{elliptic} if $|p|<\cos(\theta/2)$ (therefore, $\varphi$
fixes two points, one of them in $\mathbb{D}$ and the other in the exterior of $\overline{\DD}$). See \cite[Chapter 0]{Sh}, for instance.

\subsubsection*{Notation.} Assume $\varphi$ is a self-map of the unit disc $\DD$. In what follows, $\varphi_n$ will denote the $n$-th iterate of the map $\varphi$, that is,
$$
\varphi_n=\varphi \circ \varphi \circ \cdots \varphi \qquad  (n
\mbox{ times}),
$$
for any $n\geq 0$, where $\varphi_0$ is the identity function. It is clear that $C_{\varphi}^n=C_{\varphi_n}$ for any $n\geq 0$. If $W_{h,\varphi}$ is a bounded weighted composition operator in $\D$, it is rather straightforward that
$$
W_{h,\varphi}^n f(z)=h(z)\cdots h(\varphi_{n-1}(z))\, f(\varphi_n(z))
$$
for any $f\in \D$ and $z\in \DD$. Following \cite{HLNS}, we will denote
$$
h_{(n)}=\prod_{k=0}^{n-1} h\circ \varphi_k;
$$
where $h_{(0)}=1$ for convenience.

In what follows, we restrict our attention to weighted composition operators $W_{h,\varphi}$ acting on $\D$ induced by disc automorphisms $\varphi$. By Theorem \ref{Theorem:inner}, this implies that $h$ is a multiplier of $\D$.

\subsection{Elliptic case}

In \cite[Section 4.3]{HLNS}, the authors describe the spectrum of $W_{h,\varphi}$ acting on a large class of spaces of analytic functions whenever $h$ is in the disc algebra $A(\mathbb{D})$ and $\varphi$ is an elliptic automorphism.
Our hypotheses on $h$ in the context of $\D$ is rather more general, since $W_{h,\varphi}$ is bounded if and only if $h\in \MD$ (and the spaces $\MD$ and $A(\mathbb{D})$ are not contained in each other).
Nevertheless, it is possible to take a bit further some of the ideas developed in \cite{HLNS} and show the following result in a similar way.

\begin{thm}\label{elliptic:invert}
Suppose that $\varphi$ is an elliptic automorphism of $\DD$ with fixed point $a\in \DD$ and $W_{h,\varphi}$ a weighted composition operator on $\D$. Then
\begin{enumerate}
\item either there exists a positive integer $j$ such that $\varphi_j(z)=z$ for all $z\in \DD$, in which case, if $m$ is the smallest  such integer, then
    $$\sigma(W_{h,\varphi})=\overline{\{\lambda:\; \lambda^m=h_{(m)}(z), z\in \DD\}},$$
\item or  $\varphi_n\neq {\rm Id}$ for every $n$ and, if $W_{h,\varphi}$ is invertible, then
    $$\sigma(W_{h,\varphi})=\{\lambda:\; |\lambda|=|h(a)|\}.$$
\end{enumerate}
\end{thm}

\beginpf

The proof of (1)   goes as in \cite[Theorem 4.11]{HLNS}. The only minor change concerns  the inclusion $$\sigma(W_{h,\varphi})\subset \overline{\{\lambda:\; \lambda^m=h_{(m)}(z), z\in \DD\}},$$ where a similar argument applies taking into account the fact that if $g\in \MD$ and it is bounded away from zero, then $1/g$ is also in $\MD$.
With respect to (2), just observe that $W_{h,\varphi}^{-1}=\frac{1}{h\circ \varphi^{-1}} C_{\varphi^{-1}}$ is bounded, and hence
$1/(h\circ \varphi^{-1})$ is in $\MD$ since $\varphi^{-1}$ is a disc automorphism (see Theorem \ref{Theorem:inner}). Therefore, we refer the reader to \cite{HLNS} once more.

\endpf

\subsection{Parabolic case}

Now, let us assume that $W_{h,\varphi}$ is an invertible weighted composition operator on $\D$, where $\varphi$ is a parabolic disc automorphism. The previous ideas in \cite{HLNS} to determine the spectrum of $W_{h,\varphi}$ made an extensive use of the fact that the sequence orbit $\{\varphi_n(z_0)\}$ of a point $z_0\in \DD$ is an interpolating sequence for $\mathcal{H}^{\infty}$, a space which is assumed to be contained in the multipliers of the spaces considered (see condition (C3) in \cite{HLNS}, for instance). In the case of the Dirichlet space, the interpolating sequences for the multiplier spaces $\M(\D)$ were characterized by Marshall and Sundberg \cite{MSu}, and independently by Bishop \cite{Bi}. Nevertheless, $\{\varphi_n(z_0)\}$, $z_0\in \DD$, is no longer interpolating in $\M(\D)$ and therefore, our proof completely differs from the previous ones.
Likewise, the work in \cite{Gu} made use of inner functions, which are inappropriate in the context of $\D$.

\begin{thm}\label{thm:parabolic}
Suppose that $\varphi$ is a parabolic automorphism of $\DD$ with fixed point $a\in \TT$ and $W_{h,\varphi}$ a weighted composition operator on $\D$, determined by an $h \in \MD$ that is continuous at $a$.
If $W_{h,\varphi}$ is invertible, then
\[
\sigma(W_{h,\varphi}) = \{\lambda \in \CC: |\lambda|=|h(a)|\}.
\]
\end{thm}

\beginpf
We begin by showing that the
spectrum of $W_{h,\varphi}$ is contained in the circle $ \{\lambda: |\lambda|=|h(a)|\}$.
Recall that $h \in \MD$ implies that $h \in \mathcal{H}^\infty$ and that $|h'|^2 \, dA$ is a Carleson measure for $\D$; that is, that
there exists a constant $K>0$ such that
\[
\int_\DD |h'(z)|^2 |f|^2 \, dA(z) \le K^2 \|f\|^2_\D
\]
for all $f \in \D$. We write $\|h'\|_\C$ for the least such $K$. Moreover, we see that $1/h \in \MD$ since
$W_{h,\varphi}$ is invertible, which implies that $h(a) \ne 0$.

Recalling that $(W_{h,\varphi})^n= W_{h_{(n)},\phi_n}=T_{h_{(n)}}C_{\phi_n}$, we estimate its norm.
Since the spectral radius of $C_\phi$ is 1 (see \cite{higdon}) so that for each $\epsilon>0$ we have
$\|C_{\phi_n}\| \le (1+\epsilon)^n$ for $n$ sufficiently large, it will be sufficient to consider the   operator of multiplication by
$h_{(n)}$.
For $f \in \D$ we have
\[
(h_{(n)}f)'=  h_{(n)} f' + h'_{(n)} f.
\]
The $\A^2$ norm of the first term can be estimated using the fact that for each $\epsilon>0$
there is an $m$ such that
\beq\label{eq:hnnorm}
\|h_{(n)}\|_\infty \le \|h\|_\infty^m [(1+\epsilon)|h(a)|]^{n-m}
\eeq
for all $n \ge m$, which is given in the proof of \cite[Lem.~4.2]{HLNS}.

Also
\[
\|h'_{(n)}\|_\C \le \sum_{j=0}^{n-1} \|h_{(n),j} (h \circ \phi_j)'\|_\C,
\]
where $h_{(n),j}= h_{(n)}/(h \circ \phi_j)$.

Hence, as in (\ref{eq:hnnorm}) we have for $j < m$
\[
\|h_{(n),j}\|_\infty \le \|h\|_\infty^{m-1} [(1+\epsilon)|h(a)|]^{n-m},
\]
while for $j \ge m$ we have
\[
\|h_{(n),j}\|_\infty \le \|h\|_\infty^{m} [(1+\epsilon)|h(a)|]^{n-m-1}.
\]
Also
\[
\|(h \circ \phi_j)'\|_\C = \|\phi'_j (h' \circ \phi_j)\|_\C \le   \|\phi'_j\|_\infty \|h'\|_\C.
\]
Now, we estimate $\|\phi'_j\|_\infty$ for any $j$. Let us suppose without loss of generality that $a=-1$. Hence
\[
\phi_n(z) = \frac{(2-niy)z-niy}{niyz+2+niy},
\]
where $y \in \RR\setminus \{0\}$ (see, for instance, \cite{GM} for a similar computation when $a=1$). For a M\"obius map $ \psi(z)=\frac{z-\alpha}{1-\overline\alpha z}$
we have
\[
|\psi'(z)| = \frac{1-|\alpha|^2}{|1-\overline \alpha z|^2},
\]
and in the case of $\phi_n$ we have $\alpha=1+O(1/n)$, so that $\|\phi_n'\|_\infty = O(n)$.

Putting all this together we conclude that
\[
\|h'_{(n)}\|_\C \le C n^2 [(1+\epsilon)|h(a)|]^{n},
\]
where $C$ does not depend on $n$, and hence (since we have this for all $\epsilon>0$)
\[
\limsup_{n \to \infty} \|T_{h_{(n)}}\|^{1/n} \le |h(a)|.
\]
Having bounded the spectral radius of $W_{h,\varphi}$ by $|h(a)|$, we may
similarly bound the spectral radius of its inverse by $1/|h(a)|$, using Proposition~\ref{prop:31}.
Thus the spectrum of $W_{h,\varphi}$ is contained in the circle of
radius $|h(a)|$.

We now prove that the spectrum of $W_{h,\varphi}$ is the entire circle, by showing that for each $\lambda$ with $|\lambda|=|h(a)|$ the spectral
radius   $\rho(W_{h,\varphi}-\lambda I)$ is at least $2|h(a)|$. This technique was also used in \cite{HLNS}, but there are extra complications in our case. As we already assumed, we may choose the fixed point $a$ of $\phi$ to be   $-1$.

There is a  unitary mapping $J$ from $\D=\D(\DD)$ onto $\D(\CC_+)$, the Dirichlet space of the right half-plane, induced by the conformal involution
\[
M:z \mapsto \dfrac{1-z}{1+z};
\]
an easy calculation shows that $JW_{h,\phi}J^{-1}= \widetilde W_{h \circ M,\tau}$, the weighted composition operator on $\D(\CC_+)$ induced by a parabolic
automorphism $\tau=M \circ \phi \circ M$, which takes the form $\tau: s \to s+ic$ for some $c>0$ since it fixes $\infty$.

The following is a straightforward Corollary of \cite[Thm.~3.2]{CP14}, given that the operator
$Z:=(\widetilde W_{ h \circ M,\tau}-\lambda I)^n$ is causal (which may be most easily stated by saying that if the inverse Laplace transform of a function $u$ is supported
on $(T,\infty)$ for some $T>0$, then so is the inverse Laplace transform of its image $Zu$).\\

{\em
$\bullet$ Let $\tilde h: \CC_+ \to \CC$ be holomorphic and $\tau: \CC_+ \to \CC_+$ a causal holomorphic function such that the weighted composition
$\widetilde W_{\tilde h,\tau}$ is bounded on $\D(\CC_+)$. Then for all $\lambda \in \CC$ and $n \ge 1$ the inequality
\[
\|(\widetilde W_{\tilde h,\tau}-\lambda I)^n\|_{\mathcal{H}^2(\CC_+)} \le \|(\widetilde W_{\tilde h,\tau}-\lambda I)^n\|_{\D(\CC_+)}
\]
holds, and hence $\rho(\widetilde W_{\tilde h,\tau}-\lambda I)_{\mathcal{H}^2(\CC_+)} \le \rho(\widetilde W_{\tilde h,\tau}-\lambda I)_{\D(\CC_+)}$.
}\\

Now consider the operator $\widetilde W_{h \circ M,\tau}$ on $\mathcal{H}^2(\CC_+)$, and note that it is unitarily equivalent to the
weighted composition operator $W_{\frac{1+\phi}{1+z}h,\phi}$ on $\mathcal{H}^2(\DD)$ (see, for example, \cite{KP}). We then have
\begin{eqnarray*}
\rho(W_{h,\phi}-\lambda I)_\D &=& \rho(\widetilde W_{\tilde h,\tau}-\lambda I)_{\D(\CC_+)} \\
& \ge & \rho(\widetilde W_{\tilde h,\tau}-\lambda I)_{\mathcal{H}^2(\CC_+)}  \\
& =& \rho\left( W_{\frac{1+\phi}{1+z}h,\phi} - \lambda I\right)_{\mathcal{H}^2(\DD)} \\
& \ge &  2 \left| \frac{1+\phi(z)}{1+z}\right|_{z=-1}| h(-1) | = 2 |h(-1)|,
\end{eqnarray*}
where the last assertion
is a direct consequence of an observation made in the proof of \cite[Thm.~4.3]{HLNS}
and the fact that
\[
\frac{1+\phi(z)}{1+z} \to 1 \qquad \hbox{as} \quad z \to -1,
\]
since $\phi$ is parabolic.
Since $\rho(W_{h,\phi}-\lambda I)_\D \ge 2|h(-1)|$, we have
\[
\sigma(W_{h,\varphi}) = \{\lambda \in \CC: |\lambda|=|h(-1)|\},
\]
as required.
\endpf

\subsection{Hyperbolic case}

The same method as we adopted for the parabolic case can be used to show that the
spectrum of $W_{h,\phi}$ is contained in an annulus, in the case that $\phi$ is a
hyperbolic automorphism.

\begin{thm}\label{hyperbolic}
Suppose that $\varphi$ is a hyperbolic automorphism of $\DD$ with attractive fixed point $a\in \TT$ and
repelling fixed point $b \in \TT$. Let
$W_{h,\varphi}$ be a weighted composition operator on $\D$, determined by an $h \in \MD$ that is continuous at $a$
and $b$.
If $W_{h,\varphi}$ is invertible, then
\[
\rho(W_{h,\varphi}) \le \max \{|h(a)|, |h(b)|\}/ \mu,
\]
where $\phi$ is conjugate to the automorphism
\[
\psi(z)=\frac{(1+\mu)z+(1-\mu)}{(1-\mu)z+(1+\mu)},
\]
with $0<\mu<1$. Hence
$\sigma(W_{h,\varphi})$ is contained in the annulus with radii
$\max \{|h(a)|, |h(b)|\}/ \mu$ and $\min \{|h(a)|, |h(b)|\} \mu$.
\end{thm}

\beginpf
For each $\epsilon>0$ there is an $m$ such that the estimate
\[
\|h_{(n)}\|_\infty \le \|h\|_\infty^m [(1+\epsilon)\max \{|h(a)|, |h(b)|\}]^{n-m}
\]
holds, as in \cite{HLNS}.

Similarly, for the derivative, we have
\[
\|h'_{(n)}\|_\C \le \sum_{j=0}^{n-1} \|h_{(n),j} (h \circ \phi_j)'\|_\C,
\]
where $h_{(n),j}= h_{(n)}/(h \circ \phi_j)$, and where $\|\cdot\|_\C$ has been defined in the proof of Theorem~\ref{thm:parabolic}.
Now, however,
\[
\|(h \circ \phi_j)'\|_\C = \|\phi'_j (h' \circ \phi_j)\|_\C \le   \|\phi'_j\|_\infty \|h'\|_\C,
\]
where $\|\phi'_j\|_\infty = O(\mu^{-j})$, using the fact that
\[
\psi_j(z)=\frac{(1+\mu^j)z+(1-\mu^j)}{(1-\mu^j)z+(1+\mu^j)},
\]
as in \cite[p.~38]{GM}. Also
\[
\|h_{(n),j}\|_\infty \le \|h\|_\infty^{m-1} [(1+\epsilon)\max \{|h(a)|, |h(b)|\}]^{n-m}
\]
for $j<m$ and
\[
\|h_{(n),j}\|_\infty \le \|h\|_\infty^{m} [(1+\epsilon)\max \{|h(a)|, |h(b)|\}]^{n-m-1}
\]
for $j \ge m$.

By similar arguments to those used in the proof of Theorem~\ref{thm:parabolic},
we conclude that
\[
\limsup_{n \to \infty} \|T_{h_{(n)}}\|^{1/n} \le \max \{|h(a)|, |h(b)|\}/ \mu.
\]
The final assertion follows on considering   the spectral radius of $W_{h,\varphi}^{-1}$.
\endpf

\subsection*{A final remark} Finally, regarding Theorem \ref{hyperbolic}, we would like to pose the following open question:

\noindent \emph{If $\varphi$ is a hyperbolic automorphism of $\DD$ with attractive fixed point $a\in \TT$ and
repelling fixed point $b \in \TT$ conjugated to the automorphism
$\psi(z)=\frac{(1+\mu)z+(1-\mu)}{(1-\mu)z+(1+\mu)}$,  with $0<\mu<1$ and $W_{h,\varphi}$ is a weighted composition operator on $\D$, determined by an $h \in \MD$ continuous at both $a$ and $b$, does it follow that
$$
\sigma(W_{h,\varphi})=\{z\in \mathbb{C}:\; \min \{|h(a)|, |h(b)|\} \mu \leq |z|\leq \max \{|h(a)|, |h(b)|\}/ \mu \}
$$
whenever $W_{h,\varphi}$ is invertible?}

\section*{Acknowledgements}

This work was initiated during a research visit of the first and the third
authors to the Departamento de An\'alisis Matem\'atico at Universidad
Complutense de Madrid. They are grateful for the hospitality and the
support of research grant MTM2010-16679.
The first author also acknowledges support from the London Mathematical Society (under Scheme 2).\\
In the first version of this manuscript, we proved Proposition \ref{prop:31} ourselves: 
we thank the referee for pointing out references \cite{bourdon} and \cite{zhao}, where the 
result is already proved by different methods.

\end{document}